\documentclass{amsart}

\usepackage[utf8]{inputenc}

\usepackage{url}

\usepackage[all]{xy}
\usepackage{pb-diagram, pb-xy}

\usepackage{amssymb}
\usepackage{amsthm}

\newtheorem{thm}{Theorem}[section]
\newtheorem{lem}[thm]{Lemma}
\newtheorem{cor}[thm]{Corollary}
\newtheorem{prop}[thm]{Proposition}
\newtheorem{rem}[thm]{Remark}
\newtheorem{ex}[thm]{Example}

\newtheorem{conj}[thm]{Conjecture}

\newcommand{\hooklongrightarrow}{\lhook\joinrel\longrightarrow}

\newcommand{\bbC}{{\mathbb C}}

\newcommand{\bbG}{{\mathbb G}}

\newcommand{\bbN}{{\mathbb N}}

\newcommand{\bbQ}{{\mathbb Q}}

\newcommand{\bbZ}{{\mathbb Z}}

\newcommand{\Qbar}{\overline{\mathbb Q}}

\newcommand{\calA}{{\mathcal A}}
\newcommand{\calB}{{\mathcal B}}

\newcommand{\calH}{{\mathcal H}}
\newcommand{\calJ}{{\mathcal J}}

\newcommand{\calO}{{\mathcal O}}

\newcommand{\calX}{{\mathcal X}}
\newcommand{\calY}{{\mathcal Y}}

\newcommand{\Sp}{\mathit{Sp}}
\newcommand{\SO}{\mathit{SO}}
\newcommand{\Sl}{\mathit{Sl}}
\newcommand{\Gl}{\mathit{Gl}}
\newcommand{\Rep}{{\mathit{Rep}}}
\newcommand{\Pic}{{\mathit{Pic}}}

\newcommand{\Dbc}{{{D^{\hspace{0.01em}b}_{\hspace{-0.13em}c} \hspace{-0.05em} }}}
\newcommand{\Perv}{\mathrm{Perv}}
\newcommand{\pD}{{^p \! D}}

\newcommand{\ptau}{{^{p} \! \tau}}

\newcommand{\pH}{{^p \! H}}

\newcommand{\RHom}{{R\mathcal{H}\!om}}

\newcommand{\Supp}{\mathrm{Supp}}

\newcommand{\Hom}{{\mathit{Hom}}}

\newcommand{\id}{{\mathit{id}}}

\newcommand{\one}{{\mathbf{1}}}

\newcommand{\Tnull}{{\theta_{\mathit{null}}}}
\newcommand{\Tnullg}{{\theta_{\mathit{null},g}}}

\newcommand{\overbar}[1]{\mkern 1.5mu\overline{\mkern-1.5mu#1\mkern-3mu}\mkern 3mu}

\newcommand{\bfDbar}{{\overbar{\mathbf{D}}}}
\newcommand{\bfPbar}{{\overbar{\mathbf{P}}}}

\newcommand{\etabar}{{\mkern 2.5mu \overline{\mkern -2.5mu \eta \mkern -2.5mu} \mkern 2.5mu}}
\newcommand{\sbar}{{\mkern 2mu \overline{\mkern -2mu  s } }}
\newcommand{\Sbar}{{\mkern 2mu \overline{\mkern -2mu S } }}

\newcommand{\Ybar}{{\mkern 2mu \overline{\mkern -2mu Y \mkern -2mu} \mkern 2mu}}

\newcommand{\bfD}{{\mathbf{D}}}

\newcommand{\bfN}{{\mathbf{N}}}
\newcommand{\bfP}{{\mathbf{P}}}

\begin{document}

\title[On the Tannaka group attached to the theta divisor]{On the Tannaka group attached to the Theta divisor of a generic principally polarized abelian variety}
\author{T. Kr\"amer and R. Weissauer}
\address{Centre de Math\'ematiques Laurent Schwartz\\ \'Ecole Polytechnique \\ F-91128 Palaiseau cedex, France}
\email{thomas.kraemer@polytechnique.edu}
\address{Mathematisches Institut\\ Ruprecht-Karls-Universit\"at Heidelberg\\ Im Neuenheimer Feld 288, D-69120 Heidelberg, Germany}
 \email{rweissauer@mathi.uni-heidelberg.de}

\begin{abstract}
To any closed subvariety $Y$ of a complex abelian variety one can attach a reductive algebraic group $G$ which is determined by the decomposition of the convolution powers of $Y$ via a certain Tannakian formalism. For a theta divisor $Y$ on a principally polarized abelian variety, this group $G$ provides a new invariant that naturally endows the moduli space ${\calA}_g$ of principally polarized abelian varieties of dimension~$g$ with a finite constructible stratification. We determine $G$ for a generic principally polarized abelian variety, and for $g=4$ we show that the stratification detects the locus of Jacobian varieties inside the moduli space of abelian varieties.
\end{abstract}

\maketitle

\thispagestyle{empty}

\section{Introduction}

The moduli of principally polarized abelian varieties (ppav's) play an important role in algebraic geometry. For instance, the Torelli theorem allows to recover any smooth projective curve from the theta divisor on its Jacobian~\cite{Weil} \cite{DebW}. Jacobians are very special ppav's, and the search for a good characterization of them -- the classical Schottky problem -- has a long history~\cite{BeNovikov} \cite{DeUpdate} \cite{GruSchottky}. In this context the singularities of the theta divisor play an important role that was observed already by Andreotti and Mayer~\cite{AM} \cite{CMThetaSurvey} \cite{GHThetaSurvey}. Since the study of singular varieties naturally leads to the theory of perverse sheaves, one may ask whether these also give further insight on the geometry of theta divisors. In this paper we discuss an approach along these lines that has been suggested in~\cite{WeT}.

%%%%

\medskip

\subsection{The setting} Any complex ppav~$X$ with an irreducible theta divisor $\Theta \subset X$ defines a reductive algebraic group
\[ G \;=\; G(\delta_\Theta) \] 
by a natural construction that we will explain in section~\ref{subsec:tannakian}, based on the Tannakian formalism for the convolution product of perverse sheaves from~\cite{KrWVanishing}.
The group~$G$ is a new invariant of the underlying ppav, and its finite-dimensional complex algebraic representations encode subtle geometric information. More specifically, the theta divisor corresponds to an irreducible representation $V$ of $G$, and the decomposition of the tensor powers $V^{\otimes r}$ determines how the convolution powers of the theta divisor decompose. A general result on the almost connectedness of Tannaka groups defined by perverse sheaves~\cite{WeConn} says that the connected component $G^0\subseteq G$ acts irreducibly on $V$ and that the quotient~$\pi_0 = G/G^0$ is a finite abelian group. Notice that the principal polarization determines the theta divisor only up to a translation by a point in $X(\bbC)$. The group $G$ depends on the chosen translate, so we need to make some normalization. In what follows we always choose the theta divisor~$\Theta \subset X$ to be {\em symmetric} in the sense that it is stable under the morphism~$-\id_X$. Recall that the moduli space $\calA_g$ of complex ppav's of dimension~$g$ admits a finite cover 
\[
\xymatrix@R=1.5em@M=0.5em{\calB_g \ar[d]^-{\qquad \textnormal{\normalsize of degree $2^{2g}$}} \\ \calA_g} 
% \calB_g \; \longrightarrow \; \calA_g \quad \textnormal{of degree} \quad 2^{2g}
\]
where $\calB_g$ denotes the moduli space of ppav's with a theta level structure in the sense of~\cite[sect.~I.6 and IV.7]{FC}. Each point in $\calB_g$ corresponds to a ppav in~$\calA_g$ together with a symmetric theta divisor. It follows from the general constructibility result in proposition~\ref{prop:constructible} in the appendix for $n=2$ that the open dense locus $\calB^{\hspace{0.05em} \circ}_g \subset \calB_g$ of indecomposable ppav's is a disjoint union
\[
 \calB^{\hspace{0.05em} \circ}_g \;=\; \bigsqcup_{\alpha =1}^n S_\alpha
\]
of finitely many strata such that up to isomorphism, the group $G(\delta_\Theta)=G_\alpha$ only depends on the stratum~$S_\alpha$ in which the ppav $X$ with its symmetric theta divisor~$\Theta$ lies. This stratification is constructible for the Zariski topology. Indeed, the results in the appendix imply that the groups $G_\alpha$ behave in a lower semicontinuous way in the sense that for any stratum~$S_\beta$ in the closure of~$S_\alpha$ we have a non-canonical embedding $G_\beta \hookrightarrow G_\alpha$. We also remark that the connected component of~$G(\delta_\Theta)$ does not depend on the chosen symmetric theta divisor, hence it defines a finite constructible stratification of the moduli space $\calA_g$ whose pull-back to $\calB_g$ has the previous stratification as a refinement. It is an interesting question to ask what the above stratifications look like and what Tannaka groups can occur. Not much is known in this direction. The goal of this paper is to show that for a generic ppav one has~$G(\delta_\Theta)\cong \Sp_{g!}(\bbC)$ or $G(\delta_\Theta)\cong \SO_{g!}(\bbC)$ depending on whether $g$ is even or odd. By the above semicontinuity property this also provides an upper bound on the groups for all indecomposable ppav's of dimension $g$. Furthermore, it gives a tool to check the irreducibility of perverse sheaves that arise as direct summands in convolution powers of theta divisors, and by the dimension formula in remark~\ref{rem:basicproperties}(a) below it allows to compute the Euler characteristics of these perverse sheaves.

%%%%

\subsection{Tannakian categories} \label{subsec:tannakian}

Before we discuss our results in more detail, it will be convenient to recall the Tannakian description for categories of perverse sheaves in~\cite{KrWVanishing}.  Let~$X$ be a complex algebraic variety, and let $\bfD=\Dbc(X, \bbC)$ be the derived category of bounded complexes of $\bbC$-sheaves on $X$ whose cohomology sheaves are constructible with respect to an algebraic stratification. For abelian varieties $X$ this category carries a rich structure: The group law~$a: X\times X \rightarrow X$ gives rise to a {\em convolution product}
\[
 *: \quad \bfD \times \bfD \; \longrightarrow \; \bfD, \quad K_1 * K_2 \;=\; Ra_*(K_1\boxtimes K_2)
\]
with the usual associativity and commutativity properties familiar from the tensor product of group representations; more precisely $\bfD$ is a rigid symmetric monoidal category with respect to the convolution product. The unit object of this category is the rank one skyscraper sheaf $\delta_0$ supported in the origin $0\in X(\bbC)$, and for any sheaf complex $K\in \bfD$ with Verdier dual $D(K)$ the {\em adjoint dual} in the sense of rigid symmetric monoidal categories is given by 
$$K^\vee \;=\; (-\id_X)^* D(K). $$
Note that $\bfD$ is not an abelian category but only triangulated; however, we would like to construct a {\em tensor category}, i.e.~a rigid symmetric monoidal abelian $\bbC$-linear category. Let us denote by $\bfP=\Perv(X, \bbC) \subset \bfD$ the full abelian subcategory of perverse sheaves, which is by definition the core of the middle perverse $t$-structure on the derived category~\cite[ch.~8]{HTT}. It turns out that this full abelian subcategory is almost stable under the convolution product in the following sense: Let us call a perverse sheaf {\em negligible} if its hypercohomology has Euler characteristic zero, and let $\bfN \subset \bfD$ denote the full subcategory consisting of those sheaf complexes whose perverse cohomology sheaves are all negligible. In~\cite[lemma 13.1]{KrWVanishing} it has been shown that \medskip
\begin{itemize}
\item the triangulated quotient category $\bfD \twoheadrightarrow \bfDbar = \bfD/\bfN$ exists, \medskip
\item the convolution product $*$ descends to a product $*: \bfDbar \times \bfDbar \rightarrow \bfDbar$ with respect to which $\bfDbar$ becomes a rigid symmetric monoidal category, \medskip
\item the perverse $t$-structure induces a $t$-structure on $\bfDbar$ whose core $\bfPbar \subset \bfDbar$ is a full abelian subcategory which is stable under the convolution product. \medskip
\end{itemize}

\noindent We do not lose much in passing to the quotient categories: Up to tensoring with a local system of rank one, any simple negligible perverse sheaf has the form $f^*(P)[d]$ for some $P \in \Perv(X/Y, \bbC)$, where $f: X\twoheadrightarrow X/Y$ denotes the quotient by some abelian subvariety of dimension $d=\dim(Y)>0$~\cite{WeDegenerate}~\cite[cor.~5.2]{SchnellHolonomic}. We also remark that the core $\bfPbar$ admits an alternative description as the quotient of the abelian category $\bfP$ by its Serre subcategory $\bfP \cap \bfN$ of negligible perverse sheaves. To get the desired Tannakian description, notice that for any perverse sheaf $P$ the subquotients of the convolution powers of $P\oplus P^\vee$ form an abelian subcategory $\langle P \rangle \subset \bfPbar$ which is stable under the convolution product. By \cite[cor.~13.3]{KrWVanishing}~there exists a complex algebraic group $G=G(P)$, unique up to non-canonical isomorphism, such that we have an equivalence of tensor categories
\[
 \omega: \quad \langle P \rangle \; \stackrel{\sim}{\longrightarrow} \; \Rep_\bbC(G)
\]
where the right hand side denotes the tensor category of finite-dimensional complex algebraic representations of~$G$. We have the following basic properties:

\begin{rem} \label{rem:basicproperties}
(a) The dimension of the representation defined by $Q\in \langle P \rangle$
is given by the Euler characteristic 
\medskip \nopagebreak
\[ \dim_\bbC(\omega(Q)) \;=\; \chi(Q) \;=\; \sum_{n\in \bbZ} (-1)^n \dim_\bbC(H^n(X, Q)). \]

\noindent
(b) If $P$ is semisimple, then the group $G(P)$ is reductive since it admits a semisimple faithful representation $\omega(P)$. This is in particular the case if $P=\delta_Y = i_*j_{!*} \bbC_U[d]$ is the perverse intersection cohomology sheaf supported on an irreducible closed subvariety $i:Y\hookrightarrow X$ of dimension $d$ with smooth locus $j: U\hookrightarrow Y$.

\medskip

\noindent
(c) If the subvariety $Y\subset X$ is stable under the morphism $-\id_X$, then the perverse sheaf $P=\delta_Y$ is isomorphic to its adjoint dual. Then the representation $V=\omega(P)$ is isomorphic to its dual representation $V^*=\Hom_\bbC(V, \bbC)$.
\end{rem}

\subsection{Back to theta divisors}

We are interested in the case where $Y=\Theta$ is a symmetric theta divisor on a ppav~$X$ of dimension $g>1$. For a generic ppav we will see that the choice of the symmetric theta divisor does not affect $G(\delta_\Theta)$. But before we come to this, let us emphasize that this group  behaves rather differently from classical invariants like Mumford-Tate groups. For example, for Jacobian varieties we have the following result~\cite[p.~124 and th.~14]{WeBN}, always assuming~$\Theta\subset X$ is chosen to be symmetric:

\begin{thm} \label{thm:jacobian}
For the Jacobian $X=JC$ of a smooth complex projective curve $C$ of genus $g>1$ one has
\[
  G(\delta_\Theta) \;\cong \; 
 \begin{cases}
  \; \Sp_{2g-2}(\bbC)/\epsilon_{g-1} & \textnormal{\em if $C$ is hyperelliptic}, \\
  \; \Sl_{2g-2}(\bbC)/\mu_{g-1} & \textnormal{\em otherwise},
 \end{cases}
\]
where $\mu_{g-1}$ denotes the central subgroup of $(g-1)^\mathrm{st}$ roots of unity and $\epsilon_{g-1} \subseteq \{ \pm 1\}$ 
is its intersection with the center of the symplectic group. 
\end{thm}

The proof in loc.~cit.~uses that the theta divisor is the image of the Abel-Jacobi morphism $C^{g-1} \rightarrow X$ and also shows that $\omega(\delta_\Theta)$ is the $(g-1)^\mathrm{st}$ fundamental representation of the symplectic resp.~special linear group in this case. As in~\cite{WeT} one may pose the question whether among all ppav's of dimension~$g$ only Jacobian varieties have the above Tannaka groups. More generally, one may ask for the stratification which the connected component of the Tannaka group $G(\delta_\Theta)$ defines on the moduli space $ \calA_g = \calH_g/\Sp_{2g}(\bbZ)$ of ppav's, where $\calH_g$ denotes the Siegel upper half space. Since we normalized our theta divisors to be symmetric, the general constructibility results in the appendix imply that this a finite algebraic stratification. For the open stratum the following guess seems plausible.

\begin{conj} \label{conj:group}
Every complex ppav $X$ of dimension $g>1$ with a smooth theta divisor has the Tannaka group
\[
 G(\delta_\Theta) \;\cong\; 
 \begin{cases}
  \; \SO_{g!}(\bbC) & \textnormal{\em for $g$ odd}, \\
  \; \Sp_{g!}(\bbC) & \textnormal{\em for $g$ even},
 \end{cases}
\]
and $\omega(\delta_\Theta)$ is the standard representation of this classical group.
\end{conj}

Note that for~$g=2$ this is true by theorem~\ref{thm:jacobian}, since in this case any ppav with a smooth theta divisor is the Jacobian of a smooth hyperelliptic curve. For~$g=3$ any ppav with a smooth theta divisor is the Jacobian of a smooth non-hyperelliptic curve, and once again theorem~\ref{thm:jacobian} applies (the second fundamental representation of~$\Sl_4(\bbC)/\mu_2$ is identified with the standard representation of $\SO_6(\bbC)$ under the exceptional isomorphism between these two groups). So, the first open case occurs in dimension $g=4$ where the realm of Jacobian varieties is left. 

\medskip

In what follows, we put $G(g)=\SO_{g!}(\bbC)$ if~$g$ is odd, but~$G(g)=\Sp_{g!}(\bbC)$ if~$g$ is even. Let us say that a statement holds for a {\em general ppav} of dimension~$g$ if it holds for every ppav in a Zariski-open dense subset of $\calA_g$. Thus, the theta divisor of a general ppav is smooth. The main goal of this paper is to establish the following generic version of the above conjecture.

\begin{thm} \label{thm:generic_group}
For a general complex ppav $X$ of dimension $g>1$ the theta divisor has the Tannaka group
\[ G(\delta_\Theta) \cong G(g), \]
and $\omega(\delta_\Theta)$ is the standard representation of dimension $g!$ of this group.
\end{thm}

In the first relevant case $g=4$ we have $G(g)\cong\Sp_{24}(\bbC)$, and in this case we can actually say more. Recall from~\cite{DebThetaSing} that the Andreotti-Mayer locus 
$\mathcal{N}_g \subset \calA_g$
of ppav's with a singular theta divisor is for $g\geq 4$ itself a divisor with two irreducible components. One component is the locus $\Tnullg$ of ppav's with a vanishing theta null, the other component contains the closure $\calJ_g$ of the locus of Jacobian varieties and for $g=4$ is equal to it~\cite{BeTnull}.

\begin{thm} \label{thm:g4}
For $g=4$ the locus in $\calA_4$ of all ppav's with $G(\delta_\Theta)\not \cong G(g)$ is a Zariski-closed subset whose only components of codimension one in the moduli space~$\calA_4$ are the Jacobian locus $\calJ_4$ and the theta null locus $\Tnull_{,4}$.
\end{thm}

The proof of this will be given along with our proof of theorem~\ref{thm:generic_group}, using that every divisor on $\calA_4$ intersects the locus of Jacobian varieties. Assuming theorem~\ref{thm:g4} we easily recover the following result of~\cite{KrWSchottky}, showing for $g=4$ that our Tannaka groups are closely related to the Schottky problem.

\begin{cor} \label{cor:schottky}
For $g=4$ the invariant $G=G(\delta_\Theta)$ determines both the Jacobian locus and the theta null locus 
\[
 \calJ_4 \; \subset \; \calA_4 \quad \textnormal{\em and} \quad \Tnull_{,4} \; \subset \; \calA_4.
\]
\end{cor}

{\em Proof.} We know from theorem~\ref{thm:jacobian} that a general ppav in $\calJ_4$ has the Tannaka group $G\cong \Sl_6(\bbC)/\mu_3$. On the other hand, for a general ppav with a vanishing theta null the perverse sheaf $\delta_\Theta$ has Euler characteristic $\chi(\delta_\Theta) = g! - 2 = 22$ by part (2) of proposition~\ref{prop:degen} below. 
Hence by the dimension estimates in~\cite{AEV} such a ppav cannot have the Tannaka group $\Sl_6(\bbC)/\mu_3$.
In fact,  with similar arguments as in the proof of theorem~\ref{thm:jacobian}, a degeneration into a Jacobian variety shows that a general ppav in $\Tnull_{,4}$ has the Tannaka group $G\cong \Sp_{22}(\bbC)$.
\qed

\medskip

For the proof of theorem~\ref{thm:generic_group} and~\ref{thm:g4} we proceed as follows. In section~\ref{sec:bilinear_forms} we observe that by self-duality one always has an embedding $G(\delta_\Theta) \hookrightarrow G(g)$. Thus the main task will be to find sufficiently large lower bounds on the Tannaka group of a general ppav. By constructibility (see the appendix), it suffices to do this for the generic fibre of a suitable family of ppav's. So we can use degeneration arguments, since the Tannaka group of the generic fibre contains the Tannaka groups of the perverse sheaves of nearby cycles on the special fibres~\cite[lemma~14.1]{KrWVanishing}. In section~\ref{sec:N_filtration} we gather some general properties of the monodromy filtration that allow to get hold on the nearby cycles in terms of weights.~In section~\ref{sec:degen} we apply this to three degenerations of a generic ppav where the degenerate fibre has a vanishing theta null or splits as a product of ppav's or is a Jacobian variety. In each case, the nearby cycles contain a large irreducible constituent described geometrically. This will reduce the proof of our theorem to standard arguments in the representation theory of reductive groups, see section~\ref{sec:proof_of_generic_group}. We also include in section~\ref{sec:motivic} an alternative argument which does not require the degeneration into the theta null locus but uses more detailed information about the other two degenerations; this alternative approach introduces motivic techniques that may be useful for other applications as well. 
Finally, the appendix in section~\ref{sec:constructibility} collects some general constructibility properties of the Tannaka groups defined in~\cite{KrWVanishing}.

\medskip

\section{An upper bound on the Tannaka group} \label{sec:bilinear_forms}

In this section we show that the symplectic or special orthogonal group $G(g)$ in conjecture~\ref{conj:group} is the biggest group that can possibly occur. Let $X$ be a complex ppav of dimension $g > 1$ and $\Theta\subset X$ an irreducible symmetric theta divisor. The irreducibility implies that the theta divisor is normal~\cite{EL} and that $\delta_\Theta$ is a simple perverse sheaf. Hence if
\[
 \omega: \quad \langle \delta_\Theta \rangle \; \stackrel{\sim}{\longrightarrow} \; \Rep_\bbC(G(\delta_\Theta))
 \smallskip
\]
denotes the (non-canonical) equivalence of tensor categories from section~\ref{subsec:tannakian}, then it follows that the representation $V=\omega(\delta_\Theta)$ is irreducible. Since we have chosen the theta divisor to be symmetric, we furthermore get from remark~\ref{rem:basicproperties}(c) that $V$ is isomorphic to its dual $V^*=\Hom_\bbC(V, \bbC)$. Every self-dual irreducible representation is either orthogonal or symplectic, depending on whether the trivial representation~$\one$ lies in the alternating square~$\Lambda^2(V)$ or in the symmetric square~$S^2(V)$. In the present case this leads to the following result:

\begin{lem} \label{lem:Sp_or_SO}
For any smooth symmetric theta divisor $\Theta\subset X$ on a complex ppav of dimension $g$ we have an embedding 
\[ G(\delta_\Theta) \; \hookrightarrow \; G(g) \]
such that $V=\omega(\delta_\Theta)$ is the restriction of the standard representation of $G(g)$.
\end{lem}

{\em Proof.} We have $\dim_\bbC (V) = \chi(\delta_\Theta) = g!$ by the Gauss-Bonnet formula because the theta divisor is smooth. To decide whether the representation $V$ is orthogonal or symplectic, recall that the trivial representation $\one = \omega(\delta_0)$ corresponds on the geometric side to the skyscraper sheaf $\delta_0$ of rank one supported in the origin, and that by self-duality this skyscraper sheaf enters precisely once as a direct summand in $\delta_\Theta * \delta_\Theta$. We must determine whether on this skyscraper summand, the action of the symmetry constraint 
$$ S = S_{\delta_\Theta, \delta_\Theta}: \quad \delta_\Theta *\delta_\Theta  \; \stackrel{\sim}{\longrightarrow} \; \delta_\Theta*\delta_\Theta 
$$
for the convolution product on the rigid symmetric monoidal category $\bfD$ is given by~$+1$ or by $-1$. Recall from~\cite[sect.~2.1]{WeBN} that the symmetry constraint is induced by the  one for complexes of vector spaces with the usual Koszul sign rule. Using that $\delta_\Theta = \bbC_\Theta[g-1]$, it follows by base change that $S$ acts by~$(-1)^{g-1}$ on the stalk cohomology  
$$ \calH^0(\delta_\Theta*\delta_\Theta)_0 \;\cong\; H^{2g-2}(\Theta, \bbC). $$
Since this stalk cohomology coincides with the one of the summand $\delta_0 \hookrightarrow \delta_\Theta*\delta_\Theta$, it follows that the representation $V$ is orthogonal if $g$ is odd and that it is symplectic if~$g$ is even. Hence
$$
 G(\delta_\Theta) \;\hookrightarrow\;
 \begin{cases}
 \; O_{g!}(\bbC) & \textnormal{if $g$ is odd},\\
 \; \Sp_{g!}(\bbC) & \textnormal{if $g$ is even}.
 \end{cases}
$$
This already proves the lemma in the even-dimensional case. In the odd-dimensional case it remains to show that the determinant character $\det(V)$ is trivial. If it were not, then by~\cite[prop.~10.1]{KrWVanishing} this character would correspond to a skyscraper sheaf~$\delta_x$ supported in a non-zero point $x\in X(\bbC)$. In fact $x$ would be a $2$-torsion point by self-duality. Since the symmetric theta divisor $\Theta\subset X$ is determined uniquely up to a translation by a $2$-torsion point, all even convolution powers $(\delta_\Theta)^{*2n}$ with $n\in \bbN$ are intrinsically defined. In particular, since $g!$ is an even number for $g>1$, the non-zero $2$-torsion point $x\in X(\bbC)$ with $\omega(\delta_x) = \det(V) \subset V^{\otimes g!}$ is intrinsically defined and does not depend on our choice of the symmetric theta divisor. But for monodromy reasons it is impossible to select such a point naturally on every ppav on a Zariski-open dense subset of the moduli space $\calA_g$. By the constructibility result in proposition~\ref{prop:constructible} of the appendix, it follows that the lemma holds for all ppav's in some open dense subset of~$\calA_g$. It then follows for all ppav's with a smooth theta divisor by a specialization argument, using the semicontinuity lemma~\ref{lem:semicontinuous}.
\qed

\medskip

\begin{rem} 
For an indecomposable ppav $X$ with a symmetric but not necessarily smooth theta divisor $\Theta \subset X$, the specialization argument at the end of the above proof still implies that the Tannaka group~$G(\delta_\Theta)$ is a subquotient of $G(g)$.
\end{rem}

%\medskip

\section{Local monodromy} \label{sec:N_filtration}

For convenience of the reader we collect in this section  some general facts that allow to control degenerations of perverse sheaves in terms of weights. Throughout we consider the following local algebraic setting. Let $S$ be the spectrum of a strictly Henselian discrete valuation ring (in our case it will be the strict Henselization of a smooth complex algebraic curve in the moduli space of ppav's, but with other applications in mind we include the case of mixed or positive characteristic). In the subsequent considerations we denote by~$s$ the closed point of $S$ and by~$\etabar$ a geometric point over the generic point $\eta$ of $S$. For a separated $S$-scheme of finite type
$
 f:  Y \longrightarrow  S
$
we consider the functor of nearby cycles~\cite[sect.~4]{Il}
\smallskip
\[
 \Psi: \quad
 \bfP(Y_\eta) \;=\;\Perv(Y_\eta, \Lambda) 
 \; \longrightarrow \;
 \bfP(Y_s)\;=\;\Perv(Y_s, \Lambda)
\smallskip
\]
on perverse sheaves with coefficients in $\Lambda = \Qbar_l$, where $l$ is some fixed prime which is invertible on $S$. For any $\delta \in \bfP(Y_\eta)$ the local monodromy group $G=Gal(\etabar / \eta)$ acts  naturally on $\Psi(\delta)$, see~\cite[exp.~XIII]{SGA7}. Since we are working over a strictly Henselian base, this local monodromy group is equal to the inertia group. 

\medskip

If $p\geq 0$ denotes the residue characteristic of the point $s$, 
the inertia group sits in an exact sequence
$
 1  \rightarrow  P  \rightarrow  G  \rightarrow  \prod_{l'\neq p}  \bbZ_{l'}(1)  \rightarrow  1 ,
$
where~$P$ is the wild inertia group and where $l'$ runs through the set of all primes different from~$p$. In what follows we always assume that the nearby cycles $\Psi(\delta)$ are {\em tame} in the sense that $P$ acts trivially on them, a condition which of course is void for $p=0$. The perverse sheaf~$\Psi(\delta)$ is then equipped with a natural action of the quotient~$\bbZ_l(1)$ of the tame inertia group. We denote by~$T: \Psi(\delta) \longrightarrow \Psi(\delta)$ the endomorphism induced by a topological generator $2\pi i$ of the group~$\bbZ_l(1)$. In the abelian category $\Perv(Y_s, \Lambda)$ we have as in~\cite[lemma 1.1]{Rei} a Jordan decomposition 
\medskip
\[ 
 \Psi(\delta) \;=\; \Psi_1(\delta) \oplus \Psi_{\neq 1}(\delta)
 \quad \textnormal{with} \quad
 \Psi_{\neq 1}(\delta) \;=\; \bigoplus_{\alpha \neq 1} \; \Psi_\alpha(\delta), 
\]
where for each $\alpha \in \Lambda$ the perverse subsheaf $\Psi_\alpha(\delta)$ is stable under the action of $T$ and killed by a power of $T - \alpha\cdot \id$. In what follows, we will be particularly interested in the perverse intersection cohomology sheaf $\delta = \delta_{\, Y_\eta}$.

\begin{rem} \label{rem:monodromy_thm}
Suppose the nearby cycles $\Psi(\delta_{\, Y_\eta})$ are tame. If $f: Y \rightarrow S$ is proper and if the geometric generic fibre $Y_\etabar$ is smooth, then after replacing $f$ by its base change $f': Y'\rightarrow S'$ under a finite branched cover~$S'\rightarrow S$ we can assume 
\[ H^\bullet(Y_s', \Psi_{\neq 1}(\delta_{\, Y_\eta'})) \;=\; 0. 
\]
\end{rem}

{\em Proof.} Let $S' \to S$ be the normalization of $S$ in a finite extension of the residue field of $\eta$ with generic point $\eta' \mapsto \eta$, and denote by $\Sbar \to S$ the normalization in the residue field of $\etabar$. For the base changes $Y' = Y\times_S S'$ and $\Ybar = Y\times_S \Sbar$ we then have a commutative diagram
\[
\xymatrix@R=1.5em@C=2em@M=0.3em{
\Ybar_s \ar[r] \ar@{=}[d] & \Ybar \ar[d] & Y_\etabar \ar[l] \ar[d] \\
Y'_s \ar[r] \ar@{=}[d] & Y' \ar[d] & Y'_{\eta'} \ar[l] \ar[d] \\
Y_s \ar[r] & Y & Y_\eta \ar[l]
}
\]
where the vertical identifications on the left hand side hold since we are working over a strictly Henselian base. So, as an object of $\Perv(Y_s, \Lambda) = \Perv(Y'_s, \Lambda)$, the nearby cycles $\Psi(\delta_{\, Y_\eta})$ remain unchanged if we replace our original family $Y\to S$ by the base change $Y' \to S'$. But the local monodromy operation, and hence in the tame case the Jordan decomposition, is modified under this replacement. 

\medskip

Now, for each $\alpha \in \Lambda$ some power of $T-\alpha \cdot \id$ acts trivially on~$\Psi_\alpha(\delta_{\, Y_\eta})$, and hence also on its hypercohomology. So the Jordan decomposition
\smallskip
\[
 H^\bullet(Y_s, \Psi(\delta_{\, Y_\eta})) \;=\;
 \bigoplus_\alpha \; H^\bullet(Y_s, \Psi_\alpha(\delta_{\, Y_\eta}))
\]
shows that $H^\bullet(Y_s, \Psi_{\neq 1}(\delta_{\, Y_\eta})) = 0$ iff $T$ acts unipotently on $H^\bullet(Y_s, \Psi(\delta_{\, Y_\eta}))$. The latter can be achieved after a finite branched base change~$S'\to S$ of the form considered above: Proper base change shows $ H^\bullet(Y_s, \Psi(\delta_{\, Y_\eta})) = H^\bullet(Y_\etabar, \delta_{\, Y_\etabar}) $, and since the geometric generic fibre~$Y_\etabar$ is smooth, Grothendieck's local monodromy theorem~\cite[th.~1.4]{Il} says that on these cohomology groups $T$ is quasi-unipotent. \qed

\medskip

Returning to an arbitrary perverse sheaf $\delta \in \bfP(Y_\eta)$ with tame nearby cycles, to get hold on $\Psi_1(\delta)$ we consider the nilpotent operator 
\[ 
N \;=\; \tfrac{1}{2\pi i} \log(T): \quad \Psi_1(\delta) \; \longrightarrow \; \Psi_1(\delta)(-1), 
\]
where $\tfrac{1}{2\pi i} \in \bbZ_l(-1)$ is the dual of the generator $2\pi i \in \bbZ_l(1)$. 
By~\cite[sect.~1.6]{DelW}, there exists a unique finite increasing filtration~$F_\bullet(\Psi_1(\delta))$ on the perverse sheaf~$\Psi_1(\delta)$ such that 
$ N(F_i(\Psi_1(\delta))) \subseteq  F_{i-2}(\Psi_1(\delta))(-1)$
and such that each iterate $N^i$ induces an isomorphism 
$ Gr_i(\Psi_1(\delta)) \stackrel{\sim}{\longrightarrow} Gr_{-i}(\Psi_1(\delta))(-i) 
$
of the graded pieces with respect to the filtration. In the following, we denote by $P_{-i}(\delta)$ the kernel of~$N$ on $Gr_{-i}(\Psi_1(\delta))$ for $i\geq 0$. Then 
%\medskip \nopagebreak
%
\[
 Gr_{-i}(\Psi_1(\delta)) \;\; \cong \;\; \bigoplus_{k\geq 0} \; P_{-i-2k}(-k)
\]
by loc.~cit. We represent this situation by the next diagram, where each horizontal line of the triangle contains the composition factors of the graded piece shown on the left (for the morphisms $N$ the Tate twists must be ignored).
\[
 \xymatrix@=0.1em{
 \vdots & \quad & & & & \quad \dots \\
 Gr_2(\Psi_1(\delta)) & \quad & & & { P_{-2}(\delta)(-2)} \ar[dd]^N_\cong &  \\
 Gr_1(\Psi_1(\delta)) & \quad & & { P_{-1}(\delta)(-1)} \ar[dd]^N_\cong & & \quad \dots \\
 Gr_0(\Psi_1(\delta)) & \quad & { \quad P_0(\delta) \quad } & & { P_{-2}(\delta)(-1)} \ar[dd]^N_\cong & \\
 Gr_{-1}(\Psi_1(\delta)) \;\;\; & \quad & & { \quad P_{-1}(\delta) \quad } & & \quad \dots \\
 Gr_{-2}(\Psi_1(\delta)) \;\;\; & \quad & & & { \quad P_{-2}(\delta) \quad } & \\
 \vdots & \quad & & & & \quad \dots
 }
\]
The lower boundary entries $P_0(\delta), P_ {-1}(\delta), P_{-2}(\delta), \dots$ of the triangle are the graded pieces of the {\em specialization} 
\[ 
 sp(\delta) \; = \; \ker(N: \Psi_1(\delta) \to \Psi_1(\delta)(-1)),
\]
with $P_0(\delta)$ as the top quotient. So, the graded pieces of~$sp(\delta)$ determine those of all the $Gr_i(\Psi_1(\delta))$. For mixed perverse sheaves of geometric origin~\cite[chapt.~6]{BBD} we have the following result of Gabber \cite[th~5.1.2]{BB}. 

\begin{rem} \label{rem:monodromy-and-weight}
If $\delta$ is pure of weight~$w$, each $Gr_i(\Psi_1(\delta))$ is pure of weight $w+i$ so that the monodromy filtration coincides with the weight filtration up to a shift. 
\end{rem}

Returning to the general case, let $j: Y_\etabar \rightarrow Y$ and $i: Y_s \rightarrow Y$ be the geometric generic resp.~special fibre. Recall~\cite[p.~48]{Il} that the perverse $t$-structure on $\bfD(Y)$ is defined in terms of the perverse $t$-structures on these fibres by 
\smallskip
\begin{align*}
 K \in \pD^{\leq 0}(Y) & \; \Longleftrightarrow \; i^*K \in \pD^{\leq 0}(Y_s) \;\; \textnormal{and} \;\; j^*K \in \pD^{\leq -1}(Y_\etabar), \bigskip \\ \bigskip
 K \in \pD^{\geq 0}(Y) & \; \Longleftrightarrow \; i^!K \in \pD^{\geq 0}(Y_s) \;\; \textnormal{and} \;\; j^*K \in \pD^{\geq -1}(Y_\etabar). 
\end{align*}
We denote by $\bfP(Y)$ the abelian category of perverse sheaves defined as the core of this perverse $t$-structure. By abuse of notation, for a perverse sheaf $\delta\in \bfP(Y_\eta)$ we also write~$\delta$ for the pull-back to~$\bfP(Y_\etabar)$. Then $Rj_!(\delta[1])$ and $Rj_*(\delta [1])$ are perverse and $j_{!*}(\delta[1])\in \bfP(Y)$ is by definition the image of the natural morphism between them in the abelian category of perverse sheaves, see loc.~cit.

\begin{lem} \label{lem:sp}
With notations as above,
\[ sp(\delta)  %\;=\; \pH^0(i^* Rj_*(\delta)) 
\;=\;  i^* (j_{!*} (\delta[1])) [-1]. \]  
\end{lem}

{\em Proof.} We first claim that in the triangulated category $\bfD(Y_s)=\Dbc(Y_s, \Lambda)$, the cone of the morphism $N$ is given by
\begin{equation*} \tag{{$\star$}} \label{eq:cone} 
   Cone\bigl(\Psi_1(\delta) \stackrel{N}{\longrightarrow} \Psi_1(\delta)(-1)\bigr) \;=\;  i^* Rj_*(\delta[1]).
\end{equation*}
Indeed, if we forget about weights, the cone of $N$ on $\Psi_1(\delta)$ is isomorphic to the cone of $T-1$ on $\Psi(\delta)$, because $T-1$ is an isomorphism on $\Psi_{\neq 1}(\delta)$ and on $\Psi_1(\delta)$ its kernel and cokernel are isomorphic to those of $N$. Hence~\eqref{eq:cone} follows by the same argument as in~\cite[eq.~(3.6.2) and thereafter]{Il}, using that the wild inertia group~$P$ acts trivially on the nearby cycles. Now if for $n\in \{0, 1\}$ we consider the perverse cohomology sheaves
$ sp^n(\delta) = \pH^n(i^*Rj_*(\delta))$, 
we get from~\eqref{eq:cone} an exact sequence of perverse sheaves
\[
 0 \longrightarrow sp^0(\delta) \longrightarrow \Psi_1(\delta) \stackrel{N}{\longrightarrow} \Psi_1(\delta)(-1) \longrightarrow sp^1(\delta) \longrightarrow 0
\]
which in particular implies $sp(\delta) = sp^0(\delta)$. To finish the proof, note that by the perverse truncation property of the middle perverse extension~\cite[sect.~III.5.1]{KW} we have $i^*(j_{!*}(\delta[1])) = \ptau_{<0}i^*(Rj_*(\delta[1]))$  and that by \eqref{eq:cone} the complex $i^*(Rj_*(\delta[1]))$ is concentrated in perverse cohomology degrees $-1$ and $0$. \qed

\medskip

{\em Back to the complex analytic case}. To translate the above local results to a more global setting which is closer to our applications, it will be convenient to reset our notation. For this let $S$ be a smooth complex algebraic curve and $f: Y\rightarrow S$ a morphism of complex algebraic varieties, smooth over the complement of some given point~$s\in S(\bbC)$. Consider the strict Henselization 
$
 \tilde{S} = \mathit{Spec}(\calO_{S, s}^{\, sh} \, )   \rightarrow   S
$
of~$S$ at the point $s$, and let $\eta$ be its generic point. Then, referring to the base change~$ 
\tilde{f}: \tilde{Y} = Y\times_S \tilde{S}  \rightarrow  \tilde{S} 
$
of the morphism~$f$ under this strict Henselization, we can consider the nearby cycles and specialization functors
\[ 
\xymatrix@C=2.5em@M=0.7em{
 \bfP(\tilde{Y}_\eta) \ar@<0.7ex>[r]^-\Psi \ar@<-0.7ex>[r]_-{sp} & \bfP(\tilde{Y}_s)
}
\] 
where $\tilde{Y}_\eta$ and~$\tilde{Y}_s$ denote the generic resp.~special fibres of $\tilde{f}$. Note that~$\tilde{Y}_s$ is naturally identified with the fibre $f^{-1}(s)\subset Y$. For the perverse intersection cohomology sheaf the passage between the global and the local picture is given by the following result, where $i: \tilde{Y}_s = f^{-1}(s) \hookrightarrow Y$ denotes the embedding.

\begin{cor} \label{cor:sp}
In the above situation,
\[ sp(\delta_{\, \tilde{Y}_\eta}) \;=\; i^*(\delta_{\, Y}[-1]).\]
\end{cor}

{\em Proof.} Our smoothness assumption on $f$ implies that the generic fibre~$Y_\eta$ is smooth, so that the  perverse intersection cohomology sheaf $\delta = \delta_{\, \tilde{Y}_\eta}$ is the constant sheaf up to a degree shift. 
It follows that the middle perverse extension $j_{!*}(\delta[1])$ in lemma~\ref{lem:sp} arises from the perverse intersection cohomology sheaf $\delta_{\, Y}$ on the total space $Y$ via the Henselization morphism $\tilde{Y} \rightarrow Y$, and we are done. \qed

\medskip

\section{Degenerations of theta divisors} \label{sec:degen}

To construct degenerations of theta divisors, we fix an integer $n\geq 3$ and consider the moduli space $\calA_{g,n}$ of ppav's of dimension $g$ with level $n$ structure. Recall that this moduli space is represented by a quasi-projective and smooth variety over~$\bbQ$ by~\cite[chapt.~7.3]{GIT}. Analytically it is the quotient 
\[ \calA_{g,n}(\bbC) \;=\; \calH_g / \Gamma_g(n)
\]
of the Siegel upper half space $\calH_g$ by the free action of the principal congruence subgroup $\Gamma_g(n) = \ker(\Sp_{2g}(\bbZ) \rightarrow \Sp_{2g}(\bbZ/n\bbZ))$. Let $p: \calX \rightarrow \calA_{g,n}$ be the universal abelian scheme. In what follows, we assume $n$ is divisible by a sufficiently high power of $2$. Then the zero locus of the Riemann theta function $\vartheta(\tau, z)$ on the universal covering $\calH_g \times \bbC^g$ descends to a relative divisor $\Theta \subset \calX$ by the theta transformation formula~\cite[cor.~on p.~85]{Igusa}. For any  morphism  $S\rightarrow \calA_{g,n}$ from a variety $S$ we denote by $X_S = \calX \times_{\calA_{g,n}} S \rightarrow S$ the associated principally polarized abelian scheme. By construction it carries a relative theta divisor $\Theta_S = \Theta \times_{\calA_{g,n}} S \subset X_S$.

\begin{lem} \label{lem:heat-eqn}
For any point $s\in \calA_{g,n}(\bbC)$ there exists a smooth quasi-projective complex curve $S \hookrightarrow \calA_{g,n}$ which passes through $s$ along a general tangent direction and which has the property that \smallskip
\begin{enumerate}
\item[\em (a)] the generic fibre of the family $\Theta_S \longrightarrow S$ is smooth, \smallskip
\item[\em (b)] the singular loci of the total space and of the special fibre satisfy
\[ Sing(\Theta_S)  \subseteq  Sing(\Theta_s). \]
\end{enumerate}
If the theta divisor $\Theta_s$ contains a singular point of precise multiplicity two, then the inclusion in (b) is strict for a suitable choice of $S$.
\end{lem}

{\em Proof.} Since $\calA_{g,n}$ is smooth and quasi-projective, for any point~$s\in \calA_{g,n}(\bbC)$ we can find a smooth quasi-projective curve $S\hookrightarrow \calA_{g,n}$ passing through the point $s$ in a general tangent direction. We can assume that our general curve $S$ is not contained in the locus of ppav's with singular theta divisor. So after shrinking $S$, we can assume that for all $t\in S(\bbC) \setminus \{s\}$ the theta divisor $\Theta_t$ is smooth. Then part~{\em (a)} is clear and~{\em (b)} follows from the observation that the total space $\Theta_S$ is given locally in the smooth variety $X_S$ as the zero locus of a single analytic function. 

\medskip

Explicitly, let $\Delta \subset S$ be an analytic coordinate disk with coordinate $w$ centered at $s$, and consider a local lift $h: \Delta \rightarrow \calH_g$ of $\Delta \hookrightarrow S\hookrightarrow \calA_{g,n}$. On the universal covering the divisor $\Theta_S$ is described as the locus 
\[ \bigl \{ (w, z) \in \Delta \times \bbC^g \mid F(w, z) = 0  \bigr\}  \; \subset \; \Delta \times \bbC^ g 
\]
where the analytic function $F(w, z) = \vartheta(h(w), z)$ vanishes. If a point $(w, z)$ on this locus defines a singular point of the relative theta divisor $\Theta_S$, then the gradient of~$F$ must vanish at this point. For the gradient in the variable~$z$ this implies
\[ 0 \;=\; (\nabla_z F) (w, z) \;=\; (\nabla_z \vartheta)(\tau, z)
 \quad \textnormal{for} \quad \tau \;=\; h(w). 
\]
Hence $(\tau, z)$ defines a singular point of the fibre $\Theta_t$ where $t\in S(\bbC)$ denotes the image of the point $\tau$. By our choice of the curve $S$, the only singular fibre~$\Theta_t$ is the one over the point $t=s$, hence claim {\em (b)} follows.

\medskip

This being said, if some point $(\tau, z)\in \calH_g \times \bbC^g$ defines a singular point of the theta divisor $\Theta_s$ with precise multiplicity two, then by definition  
$(\partial^2 \vartheta/\partial z_\alpha \partial z_\beta) (\tau, z) \neq 0$
for some $\alpha, \beta \in \{1,2,\dots, g\}$. Then the heat equation
\[
 \frac{\partial \vartheta}{\partial \tau_{\alpha \beta}}(\tau, z) \;=\; 2\pi i \cdot  (1+\delta_{\alpha \beta}) \cdot  \frac{\partial^2 \vartheta}{\partial z_\alpha \partial z_\beta}(\tau, z)
\]
implies for the gradient with respect to the variable $\tau$ that $(\nabla_\tau \vartheta)(\tau, z)\neq 0$. Hence, taking $S$ to be a curve which passes through the point $s$ in a sufficiently general tangent direction $u$, we obtain with notations as above that
\[
 (\partial F / \partial w) (0, z) \;=\; (u\cdot \nabla_\tau \vartheta) (\tau, z) \;\neq \; 0
\]
where $F(w, z)=\vartheta(h(w), z)$ is defined for $w$ in some coordinate disk $\Delta\subset S$ by a suitable local lift $h: \Delta \to \calH_g$ with $h(0)=\tau$. In particular, since the gradient of~$F$ does not vanish at the considered point, it follows that $(\tau, z)$ defines a smooth point of the total space $\Theta_S$ of our family. Thus, the inclusion in claim {\em (b)} is strict. \qed

\medskip

We will now use the above construction to obtain information about the Tannaka groups $G(\delta_{\Theta_\tau})$ by varying the point $\tau\in \calA_{g,n}(\bbC)$. To apply the results about the monodromy filtration in section~\ref{sec:N_filtration} we consider perverse sheaves with coefficients in~$\Lambda = \Qbar_l$ throughout, but the outcome may as well be read in the larger category of analytic perverse sheaves with coefficients in~$\Lambda = \bbC$. 

\medskip

Let $U\subset \calA_{g,n}$ be the Zariski-open dense locus of all ppav's with a smooth theta divisor. By the constructibility property in proposition~\ref{prop:constructible} of the appendix, there exist finitely many reductive groups~$G_1, \dots, G_m$ over $\Lambda$ and a stratification of~$U$ into locally closed subsets
\[
 U \;=\; \bigsqcup_{i=0}^m \; U_i 
 \;\;\; \textnormal{with} \;\;\;
 G(\delta_{\Theta_{\tau}})\cong G_i 
 \;\;\; \textnormal{for all geometric points $\tau$ in $U_i$}.
\]
Let us denote by $U_0 \subseteq U$ the open dense stratum, so $G_0 \subseteq G(g)$ is the Tannaka group of the theta divisor on a generic ppav. With the notations of lemma~\ref{lem:Sp_or_SO}, the perverse intersection cohomology sheaf of such a theta divisor corresponds to the restriction $V|_{G_0}$ of the standard representation $V$ of the group $G(g)$.

\begin{prop} \label{prop:degen}
With notations as above, the following properties hold for the Tannaka groups of smooth theta divisors.
\begin{enumerate}
\item If $g=g_1+g_2$ and if theorem~\ref{thm:generic_group} holds for ppav's of dimension $g_1$ and $g_2$, then the generic group $G_0$ has a subquotient isogenous to 
$G(g_1)\times G(g_2)$. \smallskip

\item There is a connected subgroup $H\hookrightarrow G_0$ %, a character $\chi: H \longrightarrow \Lambda^*$ of order two
and an irreducible representation~$W$ of $H$ such that 
\[
 V|_H \;=\; \begin{cases}
 W\oplus \one \oplus \one & \textnormal{\em if $g$ is even}, \\
 W \oplus \one & \textnormal{\em if $g$ is odd},
 \end{cases}
\]
where $\one$ denotes the one-dimensional trivial representation. \smallskip

\item For $g\geq 4$ and any stratum $U_i$ of codimension at most one in $U$ there is a homomorphism 
$  f: \Sl_{2g-2}(\Lambda) \longrightarrow G_i $
such that $f^*(V)$ contains the $(g-1)^\mathit{st}$ fundamental representation of the special linear group as a summand.
\end{enumerate}
\end{prop}

{\em Proof.} Consider a morphism from a smooth irreducible quasi-projective curve~$S$ to the moduli space $\calA_{g,n}$ such that the generic point of~$S$ is mapped into some stratum~$U_i$ in our stratification, and fix a geometric generic point $\etabar$ of~$S$. Passing to the strict Henselization of $S$ at a point $s\in S(\bbC)$ we can consider the perverse sheaf $\Psi(\delta_{\Theta_\etabar}) \in \Perv(X_s, \bbC)$ of nearby cycles. By~\cite[lemma~14.1]{KrWVanishing} there exists a closed embedding
\[
 G(\Psi(\delta_{\Theta_\etabar})) \; \hookrightarrow \; G(\delta_{\Theta_\etabar}) \;\cong\; G_i \; \hookrightarrow \; G(g) 
\]
such that on the Tannakian side the perverse sheaf $\Psi(\delta_{\Theta_\etabar})$ corresponds to the restriction of $V$ to the subgroup $G(\Psi(\delta_{\Theta_\etabar}))$. We will apply this as follows.

\medskip

For part (1) take $S \hookrightarrow \calA_{g,n}$ to be an embedding of a smooth curve which meets the locus of decomposable ppav's in a single point $s\in S(\bbC)$ but which is otherwise contained in the open dense stratum $U_0$ and has the properties in lemma~\ref{lem:heat-eqn}. We choose the point $s$ such that 
\[
 X_s \;=\; X_1 \times X_2 \quad \textnormal{and} \quad \Theta_s \;=\; (\Theta_1 \times X_2) \cup (X_1 \times \Theta_2),
\]
where for $\alpha \in \{1,2\}$ the $X_\alpha$ are general complex ppav's of dimension $g_\alpha$ with a symmetric theta divisor $\Theta_\alpha$. Like for any divisor with two components which intersect each other transversally along a smooth subvariety, we have an exact sequence 
$
 0  \longrightarrow  \delta_{\Theta_1 \times \Theta_2} 
  \longrightarrow  \Lambda_{\Theta_s}[g-1]
  \longrightarrow  \delta_{\Theta_1\times X_2} \oplus \delta_{X_1\times \Theta_2} 
  \longrightarrow  0 
$
of perverse sheaves. Restricting this sequence to the open dense subset~$ Y = X_s \setminus Sing(\Theta_S)$ we get a monomorphism
\[
 \delta_{\Theta_1 \times \Theta_2}|_Y \; \hookrightarrow \; \Lambda_{\Theta_s}[g-1]|_Y 
 \;=\; sp(\delta_{\Theta_\eta})|_Y
\]
of perverse sheaves, where the last equality holds by corollary~\ref{cor:sp}. 
Now by the last statement in lemma~\ref{lem:heat-eqn} we can assume that the singular locus $Sing(\Theta_S)$ is a proper closed subset of~$Sing(\Theta_s) = \Theta_1\times \Theta_2$, and in this case the open dense subset $Y$ will have non-empty intersection with $\Theta_1\times \Theta_2$. Then, via middle perverse extension it follows from the above that the semisimplification of $sp(\delta_{\Theta_\eta})$ contains a direct summand $\delta_{\Theta_1 \times \Theta_2}$. By the properties of the monodromy filtration in section~\ref{sec:N_filtration} the same then a fortiori holds for the semisimplification of the perverse sheaf of nearby cycles $\Psi(\delta_{\Theta_\etabar})$. This being said, our claim follows from the elementary observation that the Tannaka group $G(\delta_{\Theta_1 \times \Theta_2})$ is isogenous to $G(\delta_{\Theta_1}) \times G(\delta_{\Theta_2})$.

\medskip

For part (2) let $s\in \calA_{g,n}(\bbC)$ be a point corresponding to a general ppav $X_s$ with a vanishing theta null, and consider a general curve $S\hookrightarrow \calA_{g,n}$ which passes through~$s$ but is otherwise contained in $U_0$ and has the properties in lemma~\ref{lem:heat-eqn}. For the special fibre with a vanishing theta null we know from~\cite{DebThetaSing}, with the correction in~\cite[rem.~4.5]{GSMDoublePoints}, or from theorem 4.2 of loc.~cit.~that $\Theta_s$ has an isolated ordinary double point $e$ as its only singularity. So, for the Euler characteristic of the special fibre the Picard-Lefschetz formula~\cite[exp.~XV, th.~3.4]{SGA7} says
\[
 \chi(\Theta_s) \;=\; \chi(\Theta_\etabar) \;+\; (-1)^g \;=\; (-1)^{g-1} \cdot (g! - 1)
\]
because the generic fibre $\Theta_\etabar$ is smooth of dimension $g-1$ with $\chi(\Theta_\etabar) = (-1)^{g-1}\cdot g!$ by the Gauss-Bonnet theorem. The ordinary double point singularity of the special fibre is resolved by the blowup $\pi: \tilde{\Theta}_s \longrightarrow \Theta_s$ in the point $e$, with a smooth quadric as the exceptional divisor. The stalk cohomology of the direct image $R\pi_*\delta_{\tilde{\Theta}_s}$ can be computed using~\cite[exp.~XII, th.~3.3]{SGA7}, and via the decomposition theorem this leads to an exact sequence of perverse sheaves
$
 0  \longrightarrow  \kappa_g  \longrightarrow  \Lambda_{\Theta_s}[g-1]  \longrightarrow  \delta_{\Theta_s}  \longrightarrow  0
$
for the skyscraper sheaf
\[
 \kappa_g \;=\;
 \begin{cases}
  \delta_e(-\tfrac{g-2}{2}) & \textnormal{if $g$ is even}, \\
  0 & \textnormal{if $g$ is odd}.
 \end{cases}
\]
Hence 
\[
 \chi(\delta_{\Theta_s}) \;=\; \chi(\Lambda_{\Theta_s}[g-1]) - \chi(\kappa_g) \;=\;
 \begin{cases}
  g! - 2 & \textnormal{if $g$ is even}, \\
  g! - 1 & \textnormal{if $g$ is odd}.
 \end{cases}
\]
On the other hand, by remark~\ref{rem:monodromy_thm} we can assume that the nearby cycles for our degeneration coincide with the unipotent nearby cycles. Then corollary~\ref{cor:sp} implies that the semisimplification of the nearby cycles must have the form
\[
 \Psi(\delta_{\Theta_\etabar})^{ss} \;=\; \delta_{\Theta_s} \oplus \gamma
\]
for some perverse skyscraper sheaf $\gamma$ supported on the singular point $e$ of the special fibre. A look at the Euler characteristics shows that $\gamma$ has rank two if $g$ is even, and rank one if $g$ is odd. Now take a splitting 
$G_s = G(\delta_{\Theta_s}) \hookrightarrow G(\Psi(\delta_{\Theta_\etabar}))$
and let~$W$ be the irreducible representation of $G_s$ which corresponds to the simple perverse sheaf $\delta_{\Theta_s}$. In lemma~\ref{lem:theta-connected-component} below we will see, as a general fact about irreducible theta divisors, that the restriction of $W$ to the connected component $H=G_s^0 \subseteq G_s$ remains irreducible, so our claim follows.

\medskip

For part (3) consider the locus $J\subset \calA_{g,n}$ of all Jacobian varieties of smooth projective curves of genus $g$. For any given stratum $U_i \subseteq U$ of codimension at most one, the closure $\overline{U}_i \subseteq \calA_{g,n}$ has codimension at most one in $\calA_{g,n}$. Hence we know from \cite[cor.~(0.7)]{SVComponents} that the intersection $J\cap \overline{U}_i$ is non-empty and therefore of codimension at most one in $J$. We can then find a point $t\in \overline{U}_i(\bbC)$ which corresponds to the Jacobian variety of a non-hyperelliptic smooth curve, indeed the hyperelliptic locus has codimension greater than one in $J$ for $g\geq 4$. Now consider in the closure $\overline{U}_i$ a curve which meets $J$ in the chosen point $t$ and is otherwise contained in the open part $U_i\subset \overline{U}_i$. Define 
$ \varphi: S  \rightarrow \calA_{g,n} $
to be the normalization of this curve, and $s\in S(\bbC)$ a point with~$\varphi(s) = t$. Then the abelian scheme $X_S\to S$ has its geometric generic fibre $X_\etabar$ in the stratum~$U_i$ and the special fibre $X_s$ is the Jacobian of a smooth non-hyperelliptic curve. 
As in lemma~\ref{lem:heat-eqn} the singular loci of the relative theta divisor~$\Theta_S$ and of the special fibre~$\Theta_s$ satisfy 
$Sing(\Theta_S) \subseteq Sing(\Theta_s)$,
so corollary~\ref{cor:sp} implies that the specialization $sp(\delta_{\Theta_\eta})$ admits the simple perverse sheaf~$\delta_{\Theta_s}$ as a subquotient. Then the same holds for the nearby cycles $\Psi(\delta_{\Theta_\etabar})$ as well, and this proves our claim by theorem~\ref{thm:jacobian}. \qed

\medskip

\section{Proof of the main theorem} \label{sec:proof_of_generic_group}

We now explain how to deduce theorem~\ref{thm:generic_group} from proposition~\ref{prop:degen}. For a complex ppav $X$ of dimension~$g$ with a symmetric theta divisor $\Theta \subset X$ let $G=G(\delta_\Theta)$ be the corresponding Tannaka group, and consider the representation
\[
 V \;=\; \omega(\delta_\Theta) \; \in \; \Rep_\bbC(G)
\]
where $\omega: \langle \delta_\Theta \rangle \stackrel{\sim}{\longrightarrow} \Rep_\bbC(G)$ is a fibre functor as in section~\ref{subsec:tannakian}. A priori the group~$G$ does not have to be connected, but we always have the following

\begin{lem} \label{lem:theta-connected-component}
If the theta divisor $\Theta$ is irreducible, then the restriction $V|_{G^0}$ is an irreducible representation of the connected component $G^0 \subseteq G$.
\end{lem}

{\em Proof.} The irreducibility of the theta divisor implies that $\delta_\Theta$ is a simple perverse sheaf, and hence $V$ is an irreducible representation of the Tannaka group~$G$. But by the fundamental result of~\cite{WeConn} the group $G/G^0$ is abelian. Hence if the claim of our lemma were not true, then a look at the invariants under $G^0$ in $V\otimes V^\vee$ would imply that
$V\otimes \chi \cong V$
for some non-trivial character $\chi: G\rightarrow \bbC^*$ of finite order. Now by~\cite[prop.~10.1]{KrWVanishing} any such character $\chi$ corresponds to a skyscraper sheaf $\delta_x$ of rank one, supported in some non-zero torsion point $x\in X(\bbC)$. Hence the above isomorphism would geometrically correspond to an isomorphism
$\delta_\Theta * \delta_x \cong \delta_\Theta$,
meaning that $t_x(\Theta) = \Theta$ for the translation $t_x: X\rightarrow X, y\mapsto y+x$. But this is impossible for $x\neq 0$, indeed the morphism
$
 X \rightarrow  \Pic^0(X),  x  \mapsto  \calO_X( \Theta - t_x(\Theta) )
$
is an isomorphism by the definition of a principal polarization. \qed

\medskip

Now suppose that $X$ is a general ppav. By lemma~\ref{lem:Sp_or_SO} the corresponding Tannaka group admits an embedding
$
 G = G(\delta_\Theta)  \hookrightarrow  G(g) 
$
and $V=\omega(\delta_\Theta)$ is the restriction of the standard representation of the classical group on the right hand side. 

\begin{lem} \label{lem:not-a-product}
For a general ppav $X$ with Tannaka group $G=G(\delta_\Theta)$, the connected component $G^0 \subseteq G$ is simple modulo its center.
\end{lem}

{\em Proof.} For any reductive algebraic group $H$ the derived group $H^0_\mathit{der}=[H^0, H^0]$ is a connected semisimple group, and we denote by
$
 \tilde{H}  \twoheadrightarrow  H^0_\mathit{der}
$
its simply connected cover. The covering group  $\tilde{H}$ is a product of simply connected covers of simple algebraic groups. By the theory of reductive groups~\cite[cor.~8.1.6]{SpringerLAG} the connected component $H^0$ is the product (with finite intersection) of its derived group and its center, so any irreducible representation of $H^0$ remains also irreducible as a representation of the covering group~$\tilde{H}$.

\medskip

Now consider the group $G=G(\delta_\Theta)$ attached to a general ppav. If the lemma were not true, then the simply connected cover of $G$ could be written as
$
 \tilde{G} = G_1 \times G_2
$
for some positive-dimensional simply connected groups $G_1$ and $G_2$. Lemma~\ref{lem:theta-connected-component} shows that $G^0$ acts irreducibly on $V=\omega(\delta_\Theta)$, hence by the above the same holds for the covering group $\tilde{G} = G_1\times G_2$. Hence
\[
 V|_{\tilde{G}} \; = \; V_1 \otimes V_2
\]
for certain irreducible representations $V_i \in \Rep_\bbC(G_i)$. Note that since $V$ is a faithful representation of $G$, it follows from the definition of the simply connected covering group that the representations $V_i$ are non-trivial and $\dim(V_i)>1$. Now let $H\hookrightarrow G$ be a connected subgroup as in part (2) of proposition~\ref{prop:degen} such that
\[
 V|_H \; = \; W\oplus 
 \begin{cases}
  \one \oplus \one & \textnormal{if $g$ is even}, \\
  \one & \textnormal{if $g$ is odd},
 \end{cases}
\]
where $W$ is an irreducible representation of $H$. Via the commutative diagram
\[
\xymatrix@R=1.5em@C=1em@M=0.5em{
 \tilde{H} \ar@{..>}[rr]^-\exists \ar@{->>}[d] && \tilde{G} \ar[rrr] \ar@{->>}[d] &&& \Gl(V_1) \times \Gl(V_2) \ar[d] \\
 H \ar@{^{(}->}[rr] && G \ar@{^{(}->}[rr] && \Gl(V) \ar@{=}[r] & \Gl(V_1\otimes V_2)
}
\]
we can consider the restrictions of $V_1$ and $V_2$ to the covering group $\tilde{H}$, and by construction we have
\[ 
 V_1|_{\tilde{H}} \otimes V_2|_{\tilde{H}} \; = \; V|_{\tilde{H}} \; = \; W|_{\tilde{H}} \oplus 
 \begin{cases}
  \one \oplus \one & \textnormal{if $g$ is even}, \\
  \one & \textnormal{if $g$ is odd},
 \end{cases}
\]
where $W|_{\tilde{H}}$ is irreducible by the general remarks above. Therefore both $V_i|_{\tilde{H}}$ are irreducible, since otherwise more than one non-trivial direct summand would occur in their tensor product. But then, since 
$\Hom_{\tilde{H}}(\one, \, V_1|_{\tilde{H}} \otimes V_2|_{\tilde{H}})  \neq  0$,
adjunction shows that $V_1|_{\tilde{H}}$ and $V_2|_{\tilde{H}}$ are dual to each other. In particular then $V_1$ and~$V_2$ have the same dimension. This is impossible, because $\dim(V)=g!$ is not the square of a natural number for any $g>1$. 

\medskip

For the case $g=4$ we can alternatively use the following argument, only requiring part (3) of proposition~\ref{prop:degen}. Here both $V_1$ and $V_2$ would have dimension at most~$12$, since their dimensions must be non-trivial divisors of $\dim(V)=24$. Now take $H \hookrightarrow G$ to be a subgroup isogenous to~$\Sl_{2g-2}(\bbC)=\Sl_6(\bbC)$ as in part (3) of the proposition. The classification in~\cite{AEV} of low-dimensional representations shows that each $V_i|_{\tilde{H}}$ can contain only trivial representations and $6$-dimensional standard representations as irreducible constituents. In particular, looking at $V_1|_{\tilde{H}} \otimes V_2|_{\tilde{H}}$ one sees that the third fundamental representation of $\Sl_6(\bbC)$ cannot occur as a constituent of $V|_H$, and this contradicts our choice of the subgroup $H$.
\qed

\medskip

To prove theorem~\ref{thm:generic_group} for the Tannaka group $G=G(\delta_\Theta)\hookrightarrow G(g)$ of a general ppav, we want to show that its connected component cannot be a proper subgroup of~$G(g)$. By lemma~\ref{lem:theta-connected-component} this connected component $G^0$ is an irreducible subgroup in the sense that for the standard representation $V \in \Rep_\bbC(G(g))$ the restriction~$V|_{G^0}$ is irreducible. So we are in the situation of the following general lemma.

\begin{lem} \label{lem:simple}
Let $H$ be an irreducible connected subgroup of $G(g)$ which is simple modulo its center. Then
$\dim (H) \leq  g!$ or $H = G(g)$.
\end{lem}

{\em Proof.} Let $V$ denote the standard representation of the classical group $G(g)$, and suppose that $\dim(H)>g! = \dim(V)$. Then the restriction $V|_H$ must be one of the low-dimensional representations listed in~\cite[table 1]{AEV}. Since $\dim (V)=g!$ is the factorial of a natural number, we obtain that one of the following two cases must occur for some natural number $r\in \bbN$.
\begin{enumerate}
 \item[$(a)$] The group $H$ is of type $A_r$, $C_r$ or $D_r$ and acts on $V$ via its standard representation.
 
 \item[$(b)$] The group $H$ is of type $A_r$ and acts on $V$ via the symmetric or via the alternating square of its standard representation.
\end{enumerate}
Case {$(a)$} can only occur for $H = G(g)$, since $G(g)$ is itself the symplectic or special orthogonal group with $V$ as its standard representation. To deal with case~{$(b)$}, note that by~\cite[th.~3.2.14]{GW} the symmetric square of the standard representation of $A_r$ is self-dual only if $r=1$, but then it has dimension~$3 \neq \dim (V)$. The alternating square of the standard representation of $A_r$ is self-dual of dimension $g!$ only in the case~$(r,g)=(3,3)$, but then $H = G(g) = \SO_6(\bbC)$ for dimension reasons. \qed

\bigskip

To complete the proof of theorem~\ref{thm:generic_group} we must show that for a general ppav~$X$ of dimension $g$ with theta divisor $\Theta$ we have
\[
 \dim(G) \;>\; g! \quad \textnormal{for the Tannaka group} \quad G = G(\delta_\Theta).
\]
This follows by induction on $g$. Indeed, by the remarks after conjecture~\ref{conj:group} we can assume $g\geq 4$. For $g = 4$ we have $\dim(G) \geq \dim(\Sl_{2g-2}(\bbC)) = 35 >  g! = 24$ by part~(3) of proposition~\ref{prop:degen}. The induction step is then provided by part (1) of the same proposition using $\dim(G(g-1)) > g!$ for~$g\geq 5$. Since for the case $g=4$ we have only used part~(3) of the proposition, theorem~\ref{thm:g4} follows as well.

\medskip

 \section{An alternative approach} \label{sec:motivic}

One of the main steps in finding the Tannaka group of a generic theta divisor has been to show that this group is simple modulo its center, see lemma~\ref{lem:not-a-product}. In this section we discuss another proof for this statement which introduces motivic techniques that may be of independent interest.
 
 \medskip
 
Let $X$ be a complex abelian variety. As in section~\ref{sec:N_filtration} we fix a prime number $l$ and consider the categories~$\bfP(X) = \Perv(X, \Lambda)$ and $\bfD(X) = \Dbc(X, \Lambda)$ with coefficients in $\Lambda = \Qbar_l$. In the terminology of~\cite[sect.~6.2.4]{BBD} any perverse sheaf~$P\in \bfP(X)$ of geometric origin is defined over a finitely generated field~$k$, and working over the algebraic closure $\bar{k}$ of $k$ we then obtain a $\Lambda$-adic representation of $Gal(\bar{k}/k)$ on the cohomology
 \[
  H^\bullet(X, \Lambda) \oplus H^\bullet(X, P).
 \]
The derived group of the connected component of the Zariski closure of the image of this Galois representation will be called the {\em motivic group} $M(P)$. One easily checks that this motivic group does not depend on the choice of $k$ and that for~$P, Q\in \bfP(X)$ of geometric origin one has an epimorphism $M(P\oplus Q)\twoheadrightarrow M(Q)$ which is compatible with the action on $H^\bullet(X, Q)$.
 
 \medskip
 
Recall from section~\ref{subsec:tannakian} that a complex $N \in \bfD(X)$ is said to be {\em negligible} if all perverse constituents of its perverse cohomology sheaves have Euler characteristic zero. A perverse sheaf without negligible subquotients is called {\em clean}. By~\cite{KrWVanishing} we have for all semisimple $P_1, P_2 \in \bfP(X)$ a decomposition $P_1 * P_2 = Q \oplus N$ into a clean perverse sheaf~$Q$ and a negligible complex $N$.

 \begin{lem} \label{lem:product}
Let $X$ be the Jacobian variety of a general curve of genus $g$ with theta divisor $\Theta$, and let $P\in \bfP(X)$ be a semisimple clean perverse sheaf of geometric origin. Then the support of the maximal clean direct summand~$Q$ in~$\delta_\Theta * P$ has dimension 
 \[
  \dim (\Supp (Q)) \; \geq \; g-1.
 \]
 The same holds if $X$ is the generic fibre of a connected family of ppav's which admits a general Jacobian variety as a special fibre.
 \end{lem}
 
 {\em Proof.} Suppose first that $X$ is the generic fibre of a connected family of ppav's which admits a general Jacobian variety $X_0$ with theta divisor $\Theta_0 \subset X_0$ as a special fibre. The nearby cycles functor $\Psi: \bfD(X) \rightarrow \bfD(X_0)$ commutes with convolution products~\cite[sect.~14]{KrWVanishing}, preserves the property of being negligible and does not increase the support dimension. Since $\delta_{\Theta_0}$ is a subquotient of $\Psi(\delta_\Theta)$, the lemma for the ppav $X$ will therefore follow from the lemma for $X_0$.
 
 \medskip
 
So it suffices to deal with the case where $X$ is the Jacobian variety of a general curve of genus $g$. We argue by contradiction. Suppose that for some negligible~$N$ we have a decomposition
 \[ 
  \delta_\Theta * P \;=\; Q \oplus N \quad \textnormal{with} \quad \dim(\Supp(Q))\leq g-2.
 \]
Since by the decomposition theorem $N$ is a direct sum of degree shifts of semisimple perverse sheaves, up to summands with vanishing hypercohomology the negligible part $N$ has the form $V^\bullet \otimes \delta_X$ for some complex $V^\bullet$ of Galois representations; this uses the classification of negligible complexes on a simple abelian variety given in~\cite[prop.~10.1]{KrWVanishing}. Since $H^\bullet(X, -)$ is a tensor functor by the K\"unneth formula, we have
\[
 H^\bullet(X, \delta_\Theta) \otimes H^\bullet(X, P) \;=\; H^\bullet(X, Q) \; \oplus \; \Bigl( V^\bullet \otimes H^\bullet(X, \delta_X)\Bigr) .
\]
Our assumption $\dim(Supp(Q)) \leq  g-2$ furthermore implies that $H^i(X, Q)=0$ for~$|i|\geq g-1$. Therefore the above equation, rewritten in terms of the intersection cohomology $IH^\bullet(\Theta, \Lambda) = H^\bullet(X, \delta_\Theta)[1-g]$ and $H^\bullet(X, \Lambda)=H^\bullet(X, \delta_X)[-g]$, shows that \medskip
\[
\bigoplus_{i\geq 0} \; IH^{i}(\Theta, \Lambda) \otimes H^{d-i-1}(X, P) \;=\; \bigoplus_{n\geq 0} \; V^{d-n} \otimes H^n(X, \Lambda) \quad \textnormal{for} \quad d\leq 1.
\]
Now $P$ being clean implies $H^{-g}(X, P)=0$, so in the last displayed equation only terms with $i<g$ occur. In these degrees we have $IH^{i}(\Theta, \Lambda) = \bigoplus_{k\geq 0} H^{i-2k}(X, \Lambda)$ by~\cite[cor.~13]{WeBN} and hence
\medskip
\[
\bigoplus_{n\geq 0} \; \bigoplus_{k\geq 0} \, H^{d-1-2k-n}(X, P) \otimes H^n(X, \Lambda) \;=\;
\bigoplus_{n\geq 0} \; V^{d-n} \otimes H^n(X, \Lambda) \quad \textnormal{for} \quad d\leq 1.
\]
Proceeding recursively by upward induction on $d$ starting at $d=\min\{i  \mid  V^i \neq 0\}$, we get
\medskip
\[
  \bigoplus_{k\geq 0} \; H^{d-1-2k}(X, P) \;=\; V^d  \ \mbox{ for } \ d\leq 1 .
\]
Now for the proof of the lemma we can replace $P$ by $P\oplus P^\vee$ and hence assume $P$ is isomorphic to its adjoint dual. Then $H^i(X, P) \cong H^{-i}(X, P)$ for all $i\in \bbZ$ implies that any irreducible representation of the motivic group $M=M(\delta_\Theta \oplus P)$ occuring in $H^\bullet(X, P)$ must occur in $V^\bullet$ as well.
 
\medskip

For a general Jacobian variety $X$ we have $M(\delta_X)=\Sp_{2g}(\Lambda)$. Furthermore the surjection $M\twoheadrightarrow M(\delta_X)$ admits a non-canonical splitting, and via this splitting we consider all the above hypercohomology groups as representations of the symplectic group.  We can assume $V^\bullet \neq 0$, because otherwise a comparison of degrees in the K\"unneth formula would lead to a contradiction. Let $\rho$ denote the half-sum of all positive roots of $Sp_{2g}(\bbC)$. Among all highest weights of the representation on $V^\bullet$ pick~$\beta$ so that the scalar product $(\rho, \beta)$ is maximal (possibly $\beta = 0$). Let $\beta_1, \dots, \beta_g$ denote the fundamental weights of $Sp_{2g}(\bbC)$. Since the highest weight $\beta_g$ occurs in~$H^g(X, \Lambda)$, we know $\beta + \beta_g$ occurs as a highest weight in 
\[ V^\bullet \otimes H^\bullet(X, \Lambda) \; \subset \; IH^\bullet(\Theta, \Lambda)[1]\otimes H^\bullet(X, P). \]
As $IH^\bullet(\Theta, \Lambda)$ only contains highest weights $\beta_i$ with $i<g$, it follows that $\beta + \beta_g$ is dominated by some weight  $\gamma + \beta_i$ where $i<g$ and where $\gamma$ is a highest weight occuring in $H^\bullet(X, P)$. Therefore
\[
 (\rho, \beta + \beta_g) \; \leq \; (\rho, \gamma + \beta_i).
\]
But by the result of the previous paragraph the weight $\gamma$ occuring in $H^\bullet(X, P)$ must occur in $V^\bullet$ as well, so $(\rho, \gamma) \leq (\rho, \beta)$ by our choice of $\beta$. Both inequalities together imply $(\rho, \beta_g) \leq (\rho, \beta_i)$. This gives the desired  contradiction since $i<g$. \qed

\medskip
 
By the above lemma we are now able to refine parts (1) and (3) of proposition~\ref{prop:degen} as follows, for~$\tilde{G}(g-1)$ denoting the universal cover of the group $G(g-1)$.

 \begin{prop} \label{prop:multiplicity}
 Let $G=G(\delta_\Theta)$ be the Tannaka group of a generic ppav $X$ of dimension $g$ with a symmetric theta divisor $\Theta \subset X$, and $V=\omega(\delta_\Theta) \in \Rep_\Lambda(G)$ the defining representation. Then there is
 \begin{enumerate}
 \item a homomorphism $ f: \Sl_{2g-2}(\Lambda) \rightarrow G$ such that $f^*(V)$ contains the $(g-1)^\mathrm{st}$ fundamental representation precisely once,
 \smallskip
 \item a homomorphism $ h: \tilde{G}(g-1) \rightarrow G$ such that $h^*(V)$ contains the standard representation precisely twice.
 \end{enumerate}
 \end{prop}

{\em Proof.} For part (1) consider a degeneration of $X$ into the Jacobian variety~$X_0$ of a general curve. Then for the degenerate theta divisor $\Theta_0 \subset X_0$, the perverse sheaf~$\delta_{\Theta_0}$ with $G(\delta_{\Theta_0})=\Sl_{2g-2}(\Lambda)$ is a subquotient of the nearby cycles $\Psi(\delta_\Theta)$, so for some complementary group~$H$ we have homomorphisms
 \[
  \Sl_{2g-2}(\Lambda) \times H \; \longrightarrow \;  G(\Psi(\delta_\Theta)) \; \hooklongrightarrow \; G(\delta_\Theta) \;=\; G
 \]
 where the first arrow is the universal covering and the second one comes from the tensor functoriality of the nearby cycles~\cite[lemma~14.1]{KrWVanishing}. Let $f$ be the restriction to~$\Sl_{2g-2}(\Lambda) \times \{1\}$ of the composite homomorphism. As in part (3) of proposition~\ref{prop:degen} the semisimplification of the nearby cycles has the form
 \[
  (\Psi(\delta_\Theta))^{ss} \;=\; \delta_{\Theta_0} \oplus R \quad \textnormal{with} \quad R\in \bfP(X_0)
  \quad \textnormal{and} \quad \Supp(R) \;\subseteq\; Sing (\Theta_0).
 \]
We know from~\cite[cor.~13]{WeBN} that $\omega(\delta_{\Theta_0}) = \beta_{g-1} \boxtimes \one \in  \Rep_\Lambda(\Sl_{2g-2}(\Lambda) \times H)$ is the~$(g-1)^\mathrm{st}$ fundamental representation of the special linear group. If claim~(1) were not true, we could find a clean direct summand $Q \subseteq R$ such that~$\omega(Q) = \beta_{g-1} \boxtimes W$ for some $W\in \Rep_\Lambda(H)$. But then 
\[ 
 \one \boxtimes W \; \hookrightarrow \; \bigl( \beta_{g-1} \boxtimes \one \bigr) \otimes \bigl( \beta_{g-1} \boxtimes W \bigr) \;=\; \omega\bigl(\delta_{\Theta_0} * Q\bigr)
\] 
also corresponds to some perverse sheaf $P\in \bfP(X_0)$, and then $\delta_{\Theta_0} * P = Q \oplus N$ with~$N$ negligible. Hence  lemma~\ref{lem:product} implies that $\dim(\Supp(Q)) \geq g-1$, which is impossible since by construction $\Supp(Q)\subseteq \Supp(R)$ is contained in the singular locus $Sing(\Theta_0)$.
 
 \medskip
 
 The proof of (2) is similar. Let $X$ degenerate into a product $X_1 \times X_2$ of a generic ppav of dimension $g-1$ with an elliptic curve. Then the degenerate theta divisor on this product has the form
 $ Y = (\Theta_1 \times X_2)\cup (X_1 \times \{0\})$.
 As in part (1) of proposition~\ref{prop:degen} we get
 \[
  (\Psi(\delta_\Theta))^{ss} \;=\; \delta_{Y} \oplus  2\cdot \delta_{\Theta_1 \times \{0\}} \oplus R
  \quad \textnormal{with} \quad 
  \Supp(R) \;\subsetneq\; \Theta_1 \times \{0\}.
 \]
Here $\delta_Y$ is negligible, and $\delta_{\Theta_1 \times \{0\}}$ enters with multiplicity two because it comes from the second step of the weight filtration of $sp(\delta_\Theta)$. If claim (2) fails, then as above one finds a direct summand $Q\subseteq R$, a perverse sheaf $P$ and a negligible complex~$N$ with $\delta_{\Theta_1 \times \{0\}} * P = Q \oplus N$. After a character twist on the elliptic curve we can assume by the vanishing theorem of~\cite{KrWVanishing} that for the projection $p: X_1 \times X_2 \twoheadrightarrow X_1$ the direct image $P_1=Rp_*(P)$ is perverse. Then we get
 \[ \delta_{\Theta_1} * P_1 = Q_1 \oplus N_1 \]
 where $Q_1 = Rp_*(Q)$ is supported on a strict subset of $\Theta_1$ and where $N_1 = Rp_*(N)$ is negligible. Hence again lemma~\ref{lem:product} leads to a contradiction. \qed
 
 \bigskip
 
 {\em Alternative proof of lemma~\ref{lem:not-a-product}, using the above results}. Let $G=G(\delta_\Theta)$ be the Tannaka group of a generic ppav of dimension $g$. Suppose that the universal covering group is a product $\tilde{G}=G_1 \times G_2$. The representation $V=\omega(\delta_\Theta)$ is then a tensor product
 \[ 
  V|_{\tilde{G}} \; \cong \; V_1 \boxtimes V_2 \quad \textnormal{for certain $V_i \in \Rep_\Lambda(G_i)$ with $\dim(V_i)>1$}.
 \]
 Choosing a lift of the homomorphism $h$ in proposition~\ref{prop:multiplicity}, we get a representation of the group~$H=\tilde{G}(g-1)$ on $V_1$ and $V_2$. A dimension estimate shows that $H$ must act trivially on one of these vector spaces. Indeed, for $g>4$ any \mbox{non-trivial} representation of~$H$ has dimension at least $(g-1)!$ by~\cite[table~1]{AEV}. So we may assume that~$H$ acts trivially on $V_2$. Then part (2) of proposition~\ref{prop:multiplicity} implies that $\dim(V_2) = 2$. By choosing a lift of the homomorphism~$f$ in the proposition, we get a representation of the group $F=\Sl_{2g-2}(\Lambda)$ on $V_1$ and $V_2$. Since for $g\geq 3$ this group does not have non-trivial representations of dimension two, it must act trivially on~$V_2$. Then $V|_F \cong (V_1|_F)^{\oplus 2}$ which contradicts part (1) of the proposition. \qed

\medskip

\begin{rem}
The same arguments show that in part (2) of proposition \ref{prop:multiplicity}, the restriction~$h^*(V)$ is isomorphic to the direct sum of trivial representations and two copies of the  standard representation
of $\tilde G(g-1)$.   
\end{rem}

\medskip

\section{Appendix: Constructibility} \label{sec:constructibility}

In this appendix we study how the Tannaka groups from~\cite{KrWVanishing} vary in families of perverse sheaves. We work in an algebraic setting, so as above we put~$\Lambda = \Qbar_l$ for some prime~$l$ and write $\bfP(-)=\Perv(-, \Lambda)$ and $\bfD(-)=\Dbc(-, \Lambda)$. 

\medskip

Let $S$ be an algebraic variety over an algebraically closed field $k$ of characteristic zero, and let $\calX \rightarrow S$ be an abelian scheme over $S$. If $\calY \hookrightarrow \calX$ is a closed subvariety that maps surjectively onto $S$, we can consider for any geometric point~$\sbar$ in $S$ the perverse intersection cohomology sheaf $\delta_{\, Y_\sbar} \in \bfP(X_\sbar)$ on the fibre $Y_\sbar \hookrightarrow X_\sbar$. The following example illustrates that the corresponding Tannaka groups in general do not depend on $\sbar$ in a constructible way.

\begin{ex} \label{ex:stratification}
If $S=E$ is an elliptic curve and $\calX=E\times S \to S$ is the constant family, then for the diagonal $\calY = \{ (e,e) \mid e\in E\}$ we have
\[
 G(\delta_{\, Y_\sbar}) \; \cong \;
 \begin{cases}
 \bbZ/n\bbZ & \textnormal{\em if $\sbar$ is a torsion point in $E$ of precise order $n$}, \\
 \bbG_m & \textnormal{\em if $\sbar$ is a point of infinite order in $E$}.
 \end{cases} 
\]
\end{ex}

\noindent In general we only have the following semicontinuity property.

\begin{lem} \label{lem:semicontinuous}
Let $\calY \hookrightarrow \calX$ be a closed subvariety which is smooth over $S$. Let $\eta \in S$ be a scheme-theoretic point, $s \in \overline{\{ \eta\}}$ a point in its closure, and choose geometric points $\etabar$ and $\sbar$ above them. Then we have an embedding 
$$G(\delta_{\, Y_\sbar})  \hookrightarrow  G(\delta_{\, Y_\etabar}).$$
\end{lem}

{\em Proof.} After base change to the reduced closed subscheme $\overline{\{ \eta\}} \hookrightarrow S$ we can assume that $\eta$ is the generic point of $S$. One then easily reduces our claim to the local situation of~\cite[lemma~14.1]{KrWVanishing}, and for the nearby cycles we have $\Psi(\delta_{Y_\etabar}) = \delta_{Y_\sbar}$ because by assumption the morphism $\calY \to S$ is smooth. \qed

\bigskip

To get a constructibility statement for the stratifications defined by the Tannaka groups, we need to impose some finiteness conditions. To simplify the notation, let us temporarily assume that $\calX=X$ is an abelian variety over $S=\mathit{Spec}(k)$.
For every~$P\in \bfP(X)$ we have a fibre functor
\[
 \omega: \quad \langle P \rangle \; \stackrel{\sim}{\longrightarrow} \; \Rep_\Lambda(G(P)).
\]
Like any character of the Tannaka group, the determinant character $\det(\omega(P))$ corresponds by \cite[prop.~10.1]{KrWVanishing} to a skyscraper sheaf $\delta_{x} \in \langle P \rangle$ supported on some point $x \in X(k)$. By abuse of notation, in what follows we write $\det(P)^n=\one$ to indicate that the determinant character has order dividing $n$, or equivalently that the above point $x$ is an $n$-torsion point. In many  applications the perverse sheaf~$P$ will be isomorphic to its adjoint dual $P^\vee$ and then $\det(P)^2 = \one$; by remark~\ref{rem:basicproperties}(c) this in particular holds if $P=\delta_Y$ is the perverse intersection cohomology sheaf supported on a symmetric closed subvariety $Y\hookrightarrow X$. 

\begin{lem} \label{lem:finiteness}
For any fixed $d, n\in \bbN$, the general linear group $\Gl_d(\Lambda)$ admits only finitely many conjugacy classes of subgroups isomorphic to Tannaka groups $G(P)$ of simple perverse sheaves $P\in \bfP(X)$ with $\chi(P)=d$ and $\det(P)^n=\one$.
\end{lem}

{\em Proof.} Let $G=G(P)$ be the Tannaka group of a simple perverse sheaf $P$ as in the lemma, and denote by $H=G^0$ the connected component of this Tannaka group. By~\cite{WeConn} there exists a finite subgroup $K\subset X$ of torsion points on $X$ such that we have
\[
 H\;=\; G\bigl(Rf_*(P) \bigr) \quad \textnormal{for the isogeny} \quad f: \; X \; \rightarrow \; X/K,
\]
and the group of connected components $\pi = G/H$ is a finite abelian group whose characters correspond to the rank one skyscraper sheaves supported in the points of~$K(\bbC)$. The restriction and induction functors for the subgroup $H\hookrightarrow G$ are given by the direct and inverse image:
\[
\xymatrix@M=0.5em{
 \Rep_\Lambda(G) \ar@<0.3pc>[r]^-{(-)|_H}  &  \Rep_\Lambda(H) \ar@<0.3pc>[l]^-{\mathit{Ind}^G_H}
}
\quad \textnormal{correspond to} \quad
\xymatrix@M=0.5em{
 \langle P \rangle \ar@<0.3pc>[r]^-{Rf_*} & \langle Rf_*(P) \rangle, \ar@<0.3pc>[l]^-{f^*}
}
\]
and looking at the restriction of the induction of a character $\chi: H \rightarrow \Lambda^*$ one sees that the adjoint action of $G$ fixes any such character. In other words, the adjoint action of $G$ induces the trivial action on the abelianization $H^{ab} = H/[H, H]$. So we may apply the following general fact:

\medskip

Let $V$ be a finite-dimensional vector space over $\Lambda$. Then for any fixed $n\in \bbN$, there exist only finitely many conjugacy classes of subgroups $G\hookrightarrow \Gl(V)$ with the property that
\begin{enumerate}
\item the restriction $V|_G$ is irreducible, 
\item the quotient $\pi = G/H$ by $H=G^0$ is a finite abelian group,
\item the adjoint action of $\pi$ on $H^{ab} = H/[H, H]$ is trivial,
\item the determinant $\det(V|_G)$ is a character of order at most $n$.
\end{enumerate}
For the proof of this general fact one may fix the subgroup $H\hookrightarrow \Gl(V)$ since the group $\Gl(V)$ contains only finitely many connected reductive subgroups up to conjugation. Using Mackey theory and the classification of non-connected reductive groups, one then shows that there are only finitely many subgroups $G\hookrightarrow \Gl(V)$ with given connected component $G^0 = H$ such that (1) -- (4) hold; for details we refer to~\cite[appendix B]{KraemerDiss}.  \qed

\medskip

Let us now return to an abelian scheme $p: \calX\rightarrow S$ whose base scheme is any variety $S$ over $k$. For $K\in \bfD(X)$ and geometric points $\sbar$ of $S$ we write $K_\sbar = i_\sbar^* (K)$ for the pull-back to the geometric fibre $i_\sbar: X_\sbar \rightarrow \calX$.

\begin{prop} \label{prop:constructible}
Let $n\in \bbN$ and $P\in \bfD(\calX)$ be such that for any geometric point~$\sbar$ the pull-back~$P_\sbar$ is a simple perverse sheaf with 
$\det(P_\sbar)^n = \one$.
Then there are reductive algebraic groups~$G_1, \dots, G_m$ and an algebraic stratification into locally closed subsets
\[
 S \;=\; \bigsqcup_{i=0}^m \; S_i 
 \;\;\; \textnormal{\em such that} \;\;\;
 G(P_\sbar) \;\cong\; G_i 
 \;\;\; \textnormal{\em for all geometric points $\sbar$ in $S_i$}.
\]
\end{prop}

\medskip

{\em Proof.} If $V$ is a finite-dimensional vector space over $\Lambda$, then every reductive algebraic subgroup of $\Gl_\Lambda(V)$ is determined uniquely by its invariants in the tensor powers $W^{\otimes r}$ of the representation $W=V\oplus V^\vee$ \cite[prop.~3.1(c)]{DelH}. Furthermore, if we only want to distinguish between finitely many given reductive subgroups up to conjugacy, then it suffices to consider only finitely many exponents $r \in \bbN$.

\medskip

In our case it follows via lemma~\ref{lem:finiteness} that the group $G(P_\sbar)$ is determined by the collection of all direct summands $\delta_0$ inside the convolution powers $(P_\sbar \oplus (P_\sbar)^\vee)^{*r}$ for finitely many $r$. To show that these direct summands depend on the geometric point~$\sbar$ in a constructible way, we may replace $S$ by an open dense subset and then proceed by induction on $\dim(S)$. So we may assume that there exists $K \in \bfD(\calX)$ such that
$$
 K_\sbar \;\cong\; P_\sbar \oplus (P_\sbar)^\vee
 \quad \textnormal{for all geometric points $\sbar$ in $S$}.
$$
Here we use the general fact~\cite[prop.~1.1.7]{KLFourier} that for any morphism $p: \calX\rightarrow S$  with smooth target $S$ and any constructible complex~$P\in \bfD(\calX)$, there is an open dense subset of $S$ such that the formation of the relative Verdier dual~$\RHom(P, p^{!}\Lambda_S)$ commutes with any base change that factors over this open subset. Now consider the relative convolution powers 
$$K(r) \;=\; Ra_{r,*}(K\boxtimes \cdots \boxtimes K)$$
defined by the $r$-fold addition morphism $a_r: \calX \times_S \cdots \times_S \calX \rightarrow \calX$ of our abelian scheme. These relative convolution powers are constructible sheaf complexes, hence in particular the dimension
$\dim_\Lambda \Hom(\delta_0, K(r)_\sbar)$
is a constructible function of the geometric point~$\sbar$ in $S$. The isomorphism
$$ (K(r))_\sbar \;\cong \;  (P_\sbar \oplus (P_\sbar)^\vee)^{*r} $$
and the remarks from the beginning of the proof thus imply that the isomorphism class of the group $G(P_\sbar)$ is a constructible function of the point $\sbar$ in $S$.
\qed

\bigskip \bigskip

\bibliographystyle{amsplain}
\bibliography{MyBibliography}

\end{document}